# Narayana Sequences for Cryptographic Applications

Krishnamurthy Kirthi

**Abstract**
This paper investigates the randomness and cryptographic properties of the Narayana series modulo p, where p is a prime number. It is shown that the period of the Narayana series modulo p is either $p^2+p+1$ (or a divisor) or $p^2-1$ (or a divisor). It is shown that the sequence has very good autocorrelation and crosscorrelation properties which can be used in cryptographic and key generation applications.

**Introduction**
There has been considerable recent interest in the Narayana sequences (e.g. [1]-[5]). Since they are closely related to Fibonacci sequences [6]-[10], one would expect many applications in data coding and cryptography, especially multiparty computation [11]-[16]. The properties of Fibonacci sequences modulo a prime have also been investigated [17]-[19]. Fibonacci sequences are also important in entropy problems of physics [20]-[25]. It is worthwhile then to determine if the use of the Narayana sequences can replace that of Fibonacci sequences in certain settings.

Narayana, who lived in the 14[th] century, proposed the following problem [1]: "A cow gives birth to a calf every year. When the calf is three years old, each calf gives birth to a calf at the beginning of each year. What is the number of progeny produced during twenty years by one cow". The sequence resulting from this problem is 1,1,1,2,3,4,6,9,13,19, ... and so on. Each number in the sequence is calculated by the summation of previous number and number three places before that in the sequence:

$$u_{n+1} = u_n + u_{n-2}$$

In this paper, the period of Narayana series modulo p, where p is a prime number, is investigated and found to be either $p^2+p+1$ (or a divisor) or $p^2-1$ (or a divisor). Some other characteristics of the Narayana sequence are presented. The investigation of the autocorrelation and crosscorrelation properties of the sequence reveals that they are good candidates for cryptographic and key generation applications.

**Narayana Series**

Narayana, an outstanding Indian mathematician of the 14th century, who was interested in summation of arithmetic series and magic squares, proved a more general summation in the middle of 14th century [3].

$$S_n^{(m)} = \frac{n(n+1)(n+2)\ldots(n+m)}{1 \cdot 2 \cdot 3 \cdot \ldots \cdot (m+1)}$$

Narayana applied the above equation to the problem of herd of cows and calves which is famous as Narayana's problem. Thus, using the above equation, he obtained [3]:

$$n = 1 + 20 + \frac{17 \cdot 18}{1 \cdot 2} + \frac{14 \cdot 15 \cdot 16}{1 \cdot 2 \cdot 3} + \frac{2 \cdot 3 \cdot 4 \cdot 5 \cdot 6 \cdot 7 \cdot 8}{1 \cdot 2 \cdot 3 \cdot 4 \cdot 5 \cdot 6 \cdot 7} = 2745$$

Narayana's problem can also be solved in the similar method that Fibonacci solved his rabbit problem. In the beginning of first year, there were two heads since one cow produced one calf. In the beginning of second and third year, the number of heads increased by one and therefore, the number of heads are 3 and 4 respectively. From the fourth year, the number of heads is defined as follows:

$$x_4 = x_3 + x_1, x_5 = x_4 + x_2, \ldots, x_n = x_{n-1} + x_{n-3},$$

since the number of cows for any year is summation of number of cows of previous year and number of calves which was born (= number of heads that were three years ago). We have the sequence

$$2,3,4,6,9,\ldots,u_{n+1} = u_n + u_{n-2}.$$

Thus, we obtain $u_{20} = 2745$ by computation. Now, we can consider the sequence

$$1,1,1,2,3,4,6,9,\ldots,u_{n+1} = u_n + u_{n-2},$$

with $n \geq 2$, $u_0 = 0$, $u_1 = 1$, $u_2 = 1$. These numbers are also called the Fibonacci Narayana numbers [3].

**Periods of Narayana Series modulo p**

Consider Narayana sequence modulo 3, the sequence is obtained as follows: 1,1,1,2,0,1,0,0,1,1,1,2,… and so on. Here, if three consecutive zeroes appear during sequence generation, the next number will be a zero. Since each of the three preceding digits can take values from 0 to p-1, the maximum period can only be $p^3$-1. Given a maximum length Narayana sequence, the digit-wise multiplication by 1 through p-1 would leave the sequence unchanged. Therefore, we obtain the result that the maximum period is $(p^3-1)/(p-1) = p^2+p+1$. When we consider only the preceding two digits, a similar argument would establish that the period can be $p^2-1$. Thus we have our central result:

**Theorem:**
Given Narayana sequence modulo p, where p is a prime number, the periods of the sequence will either be $p^2+p+1$ (or the divisor) or $p^2-1$ (or the divisor).

The table below provides the list of periods for first 50 prime numbers. Periods with multiples of $p^2+p+1$ (or the divisor) are assigned binary value -1 and periods with multiples of $p^2-1$ (or the divisor) are assigned binary value +1 and the resulting sequence is binary sequence B(n). Also, primes with even periods are assigned binary value 1 and primes with odd periods are assigned binary value 0 and the resulting sequence is binary sequence C(n).

Table 1: Primes from 3 to 223

| Prime numbers | Period Number | C(n) | In terms of p | B(n) | Prime numbers | Period Number | C(n) | In terms of p |
|---|---|---|---|---|---|---|---|---|
| 3 | 8 | 1 | p*p-1 | 1 | 97 | 3169 | 0 | **(p*p+p+1)/3** |
| 5 | 31 | 0 | **p*p+p+1** | -1 | 101 | 10303 | 0 | **p*p+p+1** |
| 7 | 57 | 0 | **p*p+p+1** | -1 | 103 | 10713 | 0 | **p*p+p+1** |
| 11 | 60 | 1 | (p*p-1)/2 | 1 | 107 | 11557 | 0 | **p*p+p+1** |
| 13 | 168 | 1 | p*p-1 | 1 | 109 | 11991 | 0 | **p*p+p+1** |
| 17 | 288 | 1 | p*p-1 | 1 | 113 | 991 | 0 | **(p*p+p+1)/13** |
| 19 | 381 | 0 | **p*p+p+1** | -1 | 127 | 2016 | 1 | (p*p-1)/8 |
| 23 | 528 | 1 | p*p-1 | 1 | 131 | 130 | 1 | (p*p-1)/132 |
| 29 | 840 | 1 | p*p-1 | 1 | 137 | 6256 | 1 | (p*p-1)/3 |
| 31 | 930 | 1 | p*p-1 | 1 | 139 | 1610 | 1 | (p*p-1)/12 |
| 37 | 342 | 1 | (p*p-1)/4 | 1 | 149 | 148 | 1 | (p*p-1)/148 |
| 41 | 1723 | 0 | **p*p+p+1** | -1 | 151 | 22800 | 1 | p*p-1 |
| 43 | 1848 | 1 | p*p-1 | 1 | 157 | 24807 | 0 | **p*p+p+1** |
| 47 | 46 | 1 | (p*p-1)/48 | 1 | 163 | 26733 | 0 | **p*p+p+1** |
| 53 | 468 | 1 | (p*p-1)/6 | 1 | 167 | 4648 | 1 | (p*p-1)/6 |
| 59 | 3541 | 0 | **p*p+p+1** | -1 | 173 | 172 | 1 | (p*p-1)/174 |
| 61 | 1240 | 1 | (p*p-1)/3 | 1 | 179 | 10680 | 1 | (p*p-1)/3 |
| 67 | 33 | 0 | $(p^2-1)/136$ | 1 | 181 | 32760 | 1 | p*p-1 |
| 71 | 5113 | 0 | **p*p+p+1** | -1 | 191 | 36673 | 0 | **p*p+p+1** |
| 73 | 2664 | 1 | (p*p-1)/2 | 1 | 193 | 37443 | 0 | **p*p+p+1** |
| 79 | 6240 | 1 | p*p-1 | 1 | 197 | 2156 | 1 | (p*p-1)/18 |
| 83 | 3444 | 1 | (p*p-1)/2 | 1 | 199 | 3960 | 1 | (p*p-1)/10 |
| 89 | 7920 | 1 | p*p-1 | 1 | 211 | 481 | 0 | **(p*p+p+1)/93** |

Table 2: Primes from primes 97 to 557

| Prime numbers | Period Number | C(n) | In terms of p | B(n) | Prime numbers | Period Number | C(n) | In terms of p |
|---|---|---|---|---|---|---|---|---|
| 97 | 3169 | 0 | **(p*p+p+1)/3** | -1 | 311 | 97033 | 0 | **p*p+p+1** |
| 101 | 10303 | 0 | **p*p+p+1** | -1 | 313 | 97968 | 1 | p*p-1 |
| 103 | 10713 | 0 | **p*p+p+1** | -1 | 317 | 100807 | 0 | **p*p+p+1** |
| 107 | 11557 | 0 | **p*p+p+1** | -1 | 331 | 54780 | 1 | (p*p-1)/2 |
| 109 | 11991 | 0 | **p*p+p+1** | -1 | 337 | 113568 | 1 | p*p-1 |
| 113 | 991 | 0 | **(p*p+p+1)/13** | -1 | 347 | 120408 | 1 | p*p-1 |
| 127 | 2016 | 1 | (p*p-1)/8 | 1 | 349 | 348 | 1 | (p*p-1)/350 |
| 131 | 130 | 1 | (p*p-1)/132 | 1 | 353 | 20768 | 1 | (p*p-1)/6 |
| 137 | 6256 | 1 | (p*p-1)/3 | 1 | 359 | 129241 | 0 | **p*p+p+1** |
| 139 | 1610 | 1 | (p*p-1)/12 | 1 | 367 | 67344 | 1 | (p*p-1)/2 |
| 149 | 148 | 1 | $(p^2-1)/148$ | 1 | 373 | 46501 | 0 | **(p*p+p+1)/3** |
| 151 | 22800 | 1 | p*p-1 | 1 | 379 | 378 | 1 | (p*p-1)/380 |
| 157 | 24807 | 0 | **p*p+p+1** | -1 | 383 | 73344 | 1 | (p*p-1)/2 |
| 163 | 26733 | 0 | **p*p+p+1** | -1 | 389 | 3783 | 0 | (p*p-1)/40 |
| 167 | 4648 | 1 | (p*p-1)/6 | 1 | 397 | 158007 | 0 | **p*p+p+1** |
| 173 | 172 | 1 | (p*p-1)/174 | 1 | 401 | 10720 | 1 | (p*p-1)/15 |
| 179 | 10680 | 1 | (p*p-1)/3 | 1 | 409 | 55760 | 1 | (p*p-1)/3 |
| 181 | 32760 | 1 | p*p-1 | 1 | 419 | 175981 | 0 | **p*p+p+1** |
| 191 | 36673 | 0 | **p*p+p+1** | -1 | 421 | 177663 | 0 | **p*p+p+1** |
| 193 | 37443 | 0 | **p*p+p+1** | -1 | 431 | 215 | 0 | (p*p-1)/864 |
| 197 | 2156 | 1 | (p*p-1)/18 | 1 | 433 | 187488 | 1 | p*p-1 |
| 199 | 3960 | 1 | (p*p-1)/10 | 1 | 439 | 193161 | 0 | **p*p+p+1** |
| 211 | 481 | 0 | **(p*p+p+1)/93** | -1 | 443 | 28099 | 0 | **(p*p+p+1)/7** |
| 223 | 12432 | 1 | (p*p-1)/4 | 1 | 449 | 14400 | 1 | (p*p-1)/14 |
| 227 | 226 | 1 | (p*p-1)/228 | 1 | 457 | 17404 | 1 | (p*p-1)/12 |
| 229 | 26220 | 1 | (p*p-1)/2 | 1 | 461 | 106260 | 1 | (p*p-1)/2 |
| 233 | 54523 | 0 | **p*p+p+1** | -1 | 463 | 6699 | 0 | (p*p-1)/32 |
| 239 | 8160 | 1 | (p*p-1)/7 | 1 | 467 | 218557 | 0 | **p*p+p+1** |
| 241 | 9680 | 1 | (p*p-1)/6 | 1 | 479 | 229921 | 0 | **p*p+p+1** |
| 251 | 63000 | 1 | p*p-1 | 1 | 487 | 237168 | 1 | p*p-1 |
| 257 | 66307 | 0 | **p*p+p+1** | -1 | 491 | 24108 | 1 | (p*p-1)/10 |
| 263 | 69168 | 1 | p*p-1 | 1 | 499 | 62250 | 1 | (p*p-1)/4 |
| 269 | 18090 | 1 | (p*p-1)/4 | 1 | 503 | 253513 | 0 | **p*p+p+1** |
| 271 | 73440 | 1 | p*p-1 | 1 | 509 | 259080 | 1 | (p*p-1) |
| 277 | 25576 | 1 | (p*p-1)/3 | 1 | 521 | 520 | 1 | (p*p-1)/522 |
| 281 | 79243 | 0 | **p*p+p+1** | -1 | 523 | 91176 | 1 | (p*p-1)/3 |
| 283 | 282 | 1 | (p*p-1)/284 | 1 | 541 | 293223 | 0 | **p*p+p+1** |
| 293 | 292 | 1 | (p*p-1)/294 | 1 | 547 | 99919 | 0 | **(p*p+p+1)/3** |
| 307 | 94557 | 0 | **p*p+p+1** | -1 | 557 | 38781 | 0 | (p*p-1)/8 |

As argued before, since the factors of p³-1 are p-1 and p²+p+1, the period of the Narayana sequence modulo p will either be p²+p+1 (or the divisor) or (p-1)(p+1) (or the divisor).

**Autocorrelation properties**

Autocorrelation is a measure of similarity between a sequence and time shifted replica of the sequence. Ideally, the autocorrelation function (ACF) should be impulsive i.e. peak value at zero time shift and zero values at all other time-shifts (i.e. side-lobes).

The first 20 bits of resulting binary sequence B(n) obtained from periods of the Narayana series modulo prime based on p²+p+1 (or the divisor) or p²-1 (or the divisor) are 1,-1,-1,1,1,1,-1,1,1,1,1,-1,1,1,1,-1,1,1,-1 and 1. Similarly, the first 20 bits of resulting sequence C(n) obtained from periods of the Narayana series modulo prime based on evens and odds are 1,0,0,1,1,1,0,1,1,1,0,1,1,1,0,1,0,0 and 1.

We first consider prime moduli and determine periodic autocorrelation properties of B(n) and C(n) to determine how good they are from the point of view of randomness. For convenience, the zeroes in C(n) sequence is changed to -1 so as to make off-peak autocorrelation as small as possible.

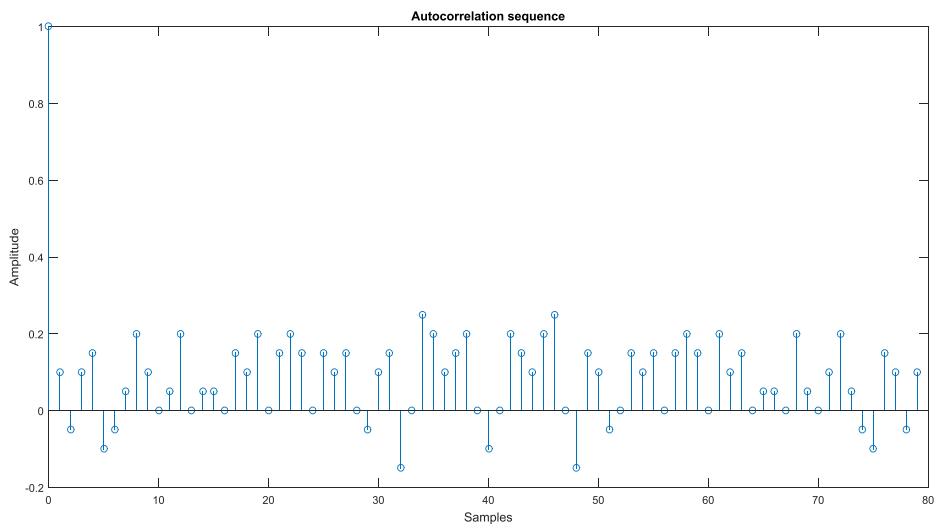

Figure 1. ACF of Binary sequence B(n) for 80 bits

The autocorrelation function is calculated using the formula:

$$ACF(k) = \frac{1}{N} \sum_{j=0}^{N-1} B_j B_{j+k}$$

where $B_j$ and $B_{j+k}$ are the binary values of sequence and time shifted version of the sequence and N is the length of sequence or period of sequence.

Figures 1 and 2 present the normalized auto correlation function of B(n) sequence for 80 and 150 bits respectively.

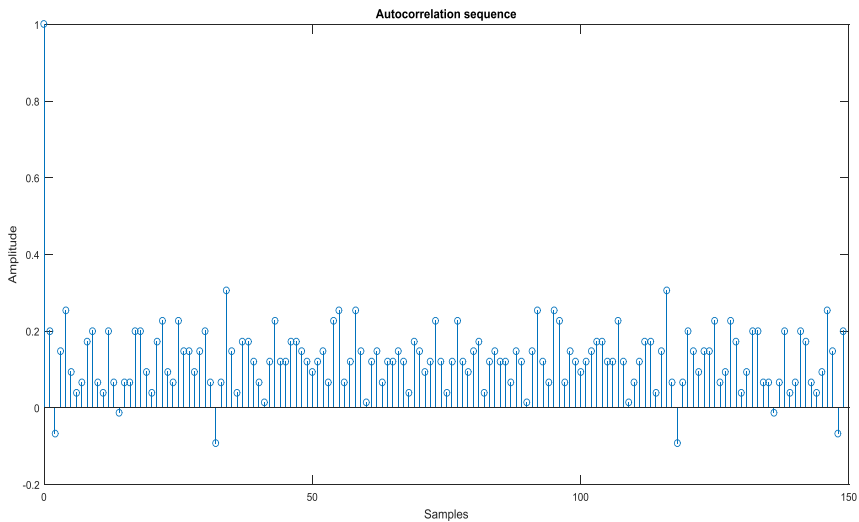

Figure 2. ACF of Binary sequence B(n) for 150 bits

Looking at Figure 1 and Figure 2, it is evident that their lobes are close to ideal autocorrelation function (ACF) with peak value at zero time shift and values close to zero at all other time-shifts (i.e. side-lobes).

Figures 3 and 4 present the normalized autocorrelation function of C(n) sequence for 100 and 140 bits respectively.

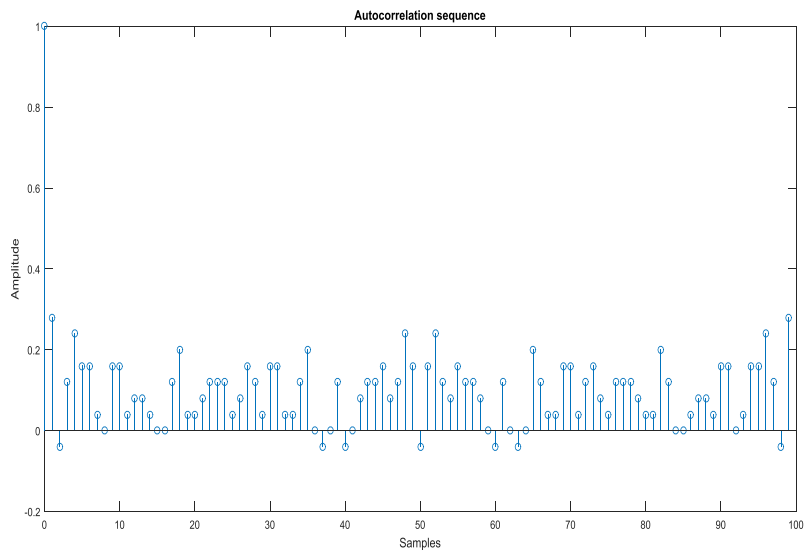

Figure 3. ACF of Binary sequence C(n) for 100 bits

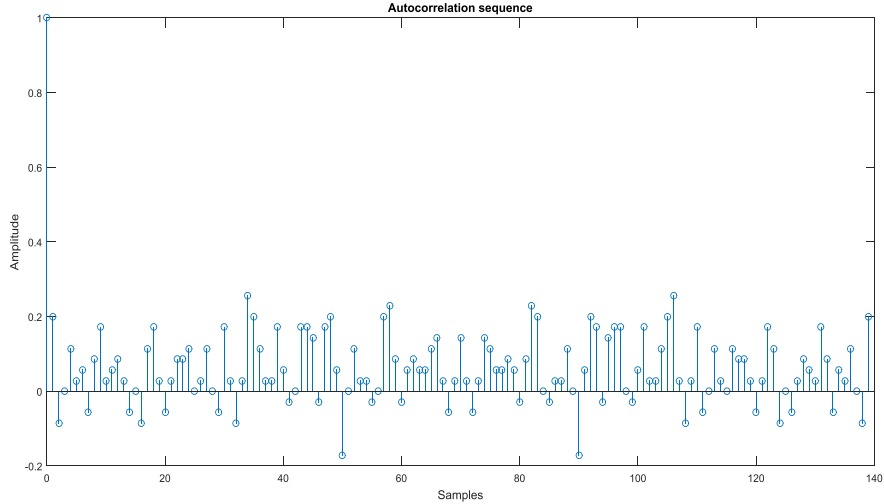

Figure 4. ACF of Binary sequence C(n) for 140 bits

Looking at Figure 3 and Figure 4, it is evident that their lobes are close to ideal autocorrelation function (ACF) with peak value at zero time shift and values close to zero at all other time-shifts (i.e. side-lobes).

Randomness may be calculated using the randomness measure, R(x), of a discrete sequence x by the expression below [17]:

$$R(x) = 1 - \frac{\sum_{k=1}^{n-1} |ACF(k)|}{N-1}$$

where *ACF(k)* is the autocorrelation function value for k and N is the period of sequence to characterize the randomness of a sequence.

According to above formula, the randomness measure of 1 indicates that the sequence is fully random whereas randomness measure of 0 indicates a constant sequence. The randomness measure for Figure 1, Figure 2, Figure 3 and Figure 4 are found to be 0.8867, 0.8662, 0.8937 and 0.9110 respectively.

**Crosscorrelation properties**
Cross orrelation is the measure of similarity between two different sequences. The cross correlation between two sequences is the complex inner product of the first sequence with a shifted version of the second sequence which indicates if the two sequences are distinct. Ideally, it is desirable to have sequences with zero crosscorrelation value at all time shifts [17]. The correlation properties of the sequences are used to detect and synchronize the communication.

Now, we consider prime moduli and determine periodic cross correlation properties of B(n) to determine how good they are from the point of view of randomness. The cross correlation function is calculated using the formula:

$$CCF(k) = \frac{1}{N} \sum_{j=0}^{N-1} A_j B_{j+k}$$

where $A_j$ and $B_{j+k}$ are the binary values of two sequences at different time intervals and N is the length of sequence or period of sequence. The peak cross correlation function value of a cross correlated sequence is calculated using the formula:

$$CCF_{peak} = \frac{1}{N} \sum_{k=1}^{N} |CCF(k)|$$

where CCF(k) is the cross correlation function value for k and N is the period of sequence. Figures 5 and 6 present the normalized cross correlation function of the B(n) sequence for 100 and 200 bits.

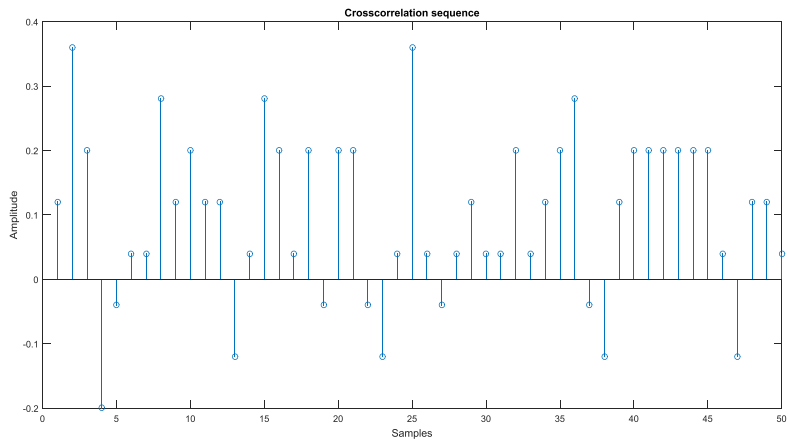

Figure 5. CCF of Binary sequence B(n) for 50 bits

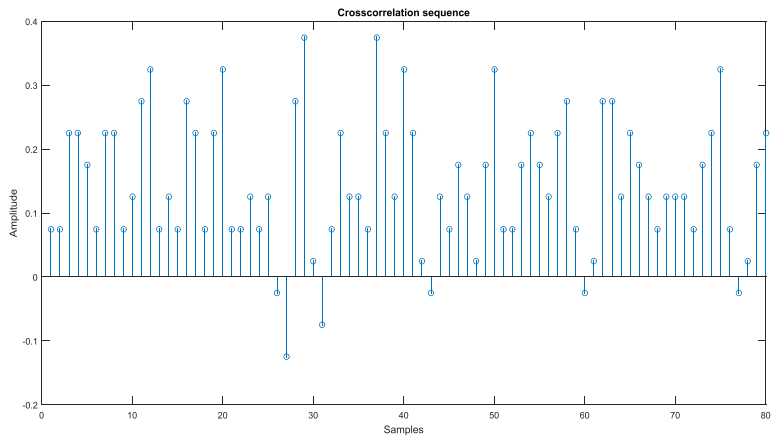

Figure 6. CCF of Binary sequence B(n) for 80 bits

Looking at Figure 5 and Figure 6, their peak crosscorrelation values are noted to be 0.375 in both cases. This value compares favorably with the peak crosscorrelation value of other pseudorandom sequences.

**Conclusion**

We have shown that the periods of Narayana sequence modulo p, p being prime number, are either $p^2+p+1$ (or a divisor) or $p^2-1$ (or a divisor). Also, we performed autocorrelation and crosscorrelation functions on Narayana series modulo prime and showed that they have good randomness properties and so they might be used for cryptographic and key distribution applications.

**References:**

1. J.-P. Allouche and T. Johnson, Narayana's cows and delayed morphisms. http://kalvos.org/johness1.html; N.J.A. Sloane, The on-line encyclopedia of integer sequences. (2008)
2. C. Flaut and V. Shpakivskyi, On generalized Fibonacci quaternions and Fibonacci-Narayana quaternions. Adv. Appl. Clifford Algebras 23: 673-688 (2013)
3. T.V. Didkivska and M.V. St'opochkina. Properties of Fibonacci-Narayana numbers. In the World of Mathematics 9 (1): 29-36 (2003) (in Ukrainian)
4. M. Bona and B.E. Sagan, On divisibility of Narayana numbers by primes. J. Integer Sequences (2005)
5. P. Barry, On a generalization of the Narayana triangle. J. Integer Sequences (2011)
6. P. Singh, The so-called Fibonacci numbers in ancient and medieval India. Historia Mathematia 12: 229-244 (1985)
7. S. Kak, The Golden Mean and the Physics of Aesthetics, arXiv: physics/0411195 (2004).
8. J. H. Thomas, Variation on the Fibonacci universal code, arXiv:cs/0701085v2 (2007)
9. M. Basu and B. Prasad, Long range variations on the Fibonacci universal code. J. of Number Theory 130: 1925-1931 (2010)
10. A. Nalli and C. Ozyilmaz, The third order variations on the Fibonacci universal code. J. of Number Theory 149: 15-32 (2015)
11. S. Kak, Oblivious transfer protocol with verification. arXiv:1504.00601 (2015)
12. S. Kak, The piggy bank cryptographic trope. Infocommunications Journal 6: 22-25 (2014)
13. L. Washbourne, A survey of P2P network security. arXiv:1504.01358 (2015)
14. S. Gangan, A review of man-in-the-middle attacks. arXiv:1504.02115 (2015)
15. S. Kak, Multiparty probability computation and verification. arXiv:1505.05081 (2015)

16. S. Kak and A. Chatterjee, On decimal sequences. IEEE Trans. on Information Theory IT-27: 647 – 652 (1981)
17. S. Gupta, P. Rockstroh, and F.E. Su, Splitting fields and periods of Fibonacci sequences modulo primes, Math. Mag. 85: 130–135 (2012)
18. M. Renault, The period, rank, and order of the (a,b)-Fibonacci sequence mod m. Math. Mag. 86: 372-380 (2013)
19. K. Kirthi, Binary GH sequences for multiparty communication, arXiv:1505.03062 (2015).
20. F. Igloi et al. Entanglement entropy of aperiodic quantum spin chains. Europhysics Letters (2007)
21. S. Kak, Quantum information and entropy. Int. Journal of Theo. Phys. 46: 860-876 (2007)
22. Y. Horibe, An entropy view of Fibonacci trees. The Fibonacci Quarterly (1982)
23. B. Varn, Optimal variable length codes: arbitrary symbol cost and equal code word probability. Information and Control 19: 289-301 (1971)
24. S. Kak, The Architecture of Knowledge. Motilal Banarsidass, Delhi (2004)
25. S. Mertens and H. Bauke, Entropy of pseudo-random number generators. Physical Review E (2004)